\documentclass[hidelinks,onefignum,onetabnum]{siamart251216}
\usepackage{amsmath,amssymb,mathtools,booktabs,xurl}
\newcommand{\C}{\mathbb C}
\newcommand{\R}{\mathbb R}

\newcommand{\Prob}{\mathbb P}
\newcommand{\fl}{\operatorname{fl}}
\newcommand{\sign}{\operatorname{sign}}
\newcommand{\gap}{\operatorname{gap}}

\newcommand{\range}{\operatorname{range}}
\newcommand{\diag}{\operatorname{diag}}
\newcommand{\INV}{\operatorname{INV}}

\newsiamthm{example}{Example}
\newsiamremark{remark}{Remark}
\headers{Logarithmic precision for Schur decomposition}{Matthew J. Colbrook}
\title{Schur decomposition with logarithmic precision\protect\\ in cubic time up to logarithmic factors}
\author{Matthew J. Colbrook\thanks{Department of Applied Mathematics and Theoretical Physics, University of Cambridge, Cambridge, United Kingdom (\email{m.colbrook@damtp.cam.ac.uk}).}}
\hypersetup{pdftitle={Schur decomposition with logarithmic precision in cubic time up to logarithmic factors},pdfauthor={Matthew J. Colbrook}}
\begin{document}
\maketitle
\begin{abstract}
How much arithmetic precision is needed to compute a Schur decomposition to a prescribed tolerance? Let $A$ be a complex matrix of order $n$ with spectral norm at most one, let $0<\delta<1/2$ be the tolerance, and let $0<p<1/2$ be the requested failure probability. Our randomised floating-point algorithm uses $O(\log(n/(\delta p)))$ mantissa bits per real component and $O(n^3\log^2(n/(\delta p)))$ arithmetic operations on every execution, with universal constants. It returns $Q$ and an exactly upper triangular $T$. With probability at least $1-p$, both $\|A-QTQ^*\|_2$ and $\|Q^*Q-I\|_2$ are at most $\delta$. We sample a small perturbation using finitely many bits and use Gaussian estimates to bound the eigenvector conditioning and eigenvalue gap. A single contour estimate bounds the effect of a rounding error through any number of remaining Newton steps. At each subdivision level, the full matrix is unitarily similar to a small perturbation of the regularised input. This allows us to bound the conditioning of each diagonal block directly from that input and to use the same precision throughout.
\end{abstract}

\section{Introduction}\label{sec:intro}

The Schur decomposition gives an orthonormal basis in which a matrix is upper triangular. The diagonal entries are the eigenvalues, and for each $k$ the first $k$ basis vectors span an invariant subspace. Many matrix computations then reduce to triangular problems. How much precision do we need to compute such a decomposition within a given operation count? The difficulty is that for a non-normal matrix even a small perturbation can move the eigenvalues or change the invariant subspaces considerably.

Let $\delta$ denote the backward error we require and $u$ the unit roundoff. At a fixed failure probability, a bound of the form
\begin{equation}\label{eq:desired-u}
 u\le (\delta/n)^C
\end{equation}
can be enforced with $O(\log(n/\delta))$ mantissa bits when the exponent $C$ is independent of dimension and tolerance. If $C$ grows like $\log n$ or $\log(1/\delta)$, the precision bound gains another logarithmic factor. Problem~3.3 in the Simons workshop report \cite{AmselEtAl2026} asks whether this increase can be avoided. We consider the explicit Schur problem, IE-08, in Townsend's problem repository \cite{TownsendRepository}.
This paper continues the project \mbox{\url{https://nla-explained.com/}}, and provides a human digestion of a proof found by AI.

Throughout, $\|\cdot\|$ denotes the Euclidean norm for vectors and the induced spectral norm for matrices, and $\|\cdot\|_F$ denotes the Frobenius norm. The matrix $X^*$ is the conjugate transpose of $X$, and $I$ is an identity matrix of the appropriate order. For a nonsingular square matrix $X$, let $\kappa_2(X)=\|X\|\|X^{-1}\|$. We denote the spectrum of a square matrix $X$ by $\sigma(X)$. For any matrix $X$, we index its singular values $\sigma_j(X)$ in nonincreasing order. Coordinate vectors are denoted by $e_i$. Logarithms without a subscript are natural.

Real arithmetic is rounded to nearest with an integer mantissa length $b\ge1$, giving unit roundoff $u=2^{-b}$. We use the relative-error model for every operation, including square roots, and assume sufficient exponent range to exclude underflow and overflow. Complex operations are evaluated through a fixed number of real operations, and all random choices use independent unbiased bits. The requests $\delta,p$ are exact control inputs for parameter selection. We use $\fl(\cdot)$ for the result of the specified kernel, including rounding at intermediate steps.

\begin{theorem}\label{thm:main}
There are a randomised floating-point algorithm and a universal constant $C>0$ with the following guarantees. For every integer $n\ge1$, every matrix $A\in\C^{n\times n}$ satisfying $\|A\|\le1$, every tolerance $0<\delta<1/2$, and every requested failure probability $0<p<1/2$, the algorithm needs at most
\[
 Cn^3\log^2(n/(\delta p))
\]
arithmetic operations and at most $C\log(n/(\delta p))$ mantissa bits per real component. The same constant applies for every $n,\delta,p$. It returns stored matrices $Q,T\in\C^{n\times n}$, with $T$ exactly upper triangular, and jointly with probability at least $1-p$,
\begin{equation}\label{eq:target}
 \|A-QTQ^*\|\le\delta,
 \qquad \|Q^*Q-I\|\le\delta.
\end{equation}
No diagonalisability or eigenvalue-separation assumption is imposed on the input. Every execution satisfies the work bound, whether or not the residual guarantee succeeds.
\end{theorem}

Both error bounds in \eqref{eq:target} are measured for the exact matrices represented by the stored entries. They include the initial rounding of the input and the errors in computing the basis and triangular factor. \Cref{thm:main} has also been formally verified in Lean \cite{ColbrookPalomar2026}.

\subsection{Perturbation, iteration and spectral subdivision}
A bound on the norm of a matrix gives no bound on the conditioning of its eigenvectors. For example,
\[
 A_\epsilon=\begin{pmatrix}0&1\\ \epsilon&0\end{pmatrix}
\]
approaches a defective matrix as $\epsilon\to0$. For real $\epsilon>0$, its eigenvectors approach parallel directions, although $\|A_\epsilon\|=1$ when $\epsilon\le1$. The theorem must allow eigenvectors that are almost parallel, as well as a defective matrix with no eigenvector basis.

We reserve part of the backward-error tolerance for a small random perturbation. The Gaussian regularisation theorem of Banks, Garza-Vargas, Kulkarni, and Srivastava \cite[Theorem~3.6]{BanksEtAl2023} gives polynomial bounds on the eigenvector condition number and reciprocal eigenvalue gap for the perturbed matrix with high probability. These bounds determine the accuracies we need later, without requiring us to compute either spectral quantity.

We divide the regularised spectrum using the matrix sign function. To compute it, Newton's iteration replaces $X$ by $(X+X^{-1})/2$. How far can a rounding error at one step move the final iterate? Multiplying separate one-step estimates can force extra logarithmic factors into the precision. We bound the whole remaining composition with a single multiplier, then use a telescoping identity to add the effects of the local errors.

At each subdivision we approximate an invariant subspace, change basis and discard the small lower-left block. We then face the same problem on the two diagonal blocks. If we pass a condition estimate from parent to child, we risk acquiring another factor at every split. We retain the full upper-block-triangular matrix in the analysis. This matrix is unitarily similar, at each level, to a small perturbation of the regularised input. Compressing its spectral projectors bounds those of the diagonal blocks, so the conditioning of each current block is controlled directly from the regularised input.

Sections~\ref{sec:sign} and \ref{sec:conditioning} give the iteration and conditioning bounds. Sections~\ref{sec:random}--\ref{sec:subspace} choose the parameters and analyse a single split. We give the algorithm in Section~\ref{sec:algorithm}, follow the errors through successive levels in Section~\ref{sec:global}, and finish the proof in Section~\ref{sec:completion}. Section~\ref{sec:numerics} measures the arithmetic work with the explicit schedule.

\subsection{Relation to earlier work}
Roberts \cite{Roberts1980} used the matrix sign function to separate the spectrum across the imaginary axis, with applications to model reduction and algebraic Riccati equations. Byers, He, and Mehrmann \cite{ByersHeMehrmann1997} studied the effect of sign-function errors on invariant-subspace computations, with particular attention to ill-conditioned Newton iterates. Bai and Demmel \cite{BaiDemmel1998} analysed perturbations of the matrix sign function and the stability of sign-based divide and conquer. For background on the matrix sign function and its spectral projectors, see Higham \cite[Chapter~5]{Higham2008}.

Banks et al.\ \cite{BanksEtAl2023} introduced pseudospectral shattering to analyse diagonalisation in nearly matrix multiplication time. Their regularisation theorem gives our starting matrix. Their analysis also uses contour integrals to bound matrix iterations. We use spectral divide and conquer with classical cubic inversion, controlling all remaining Newton steps by a single uniform bound. Combined with the diagonal-block conditioning estimate, this gives the precision in \eqref{eq:desired-u}. Banks, Garza-Vargas, and Srivastava \cite{BanksGarzaVargasSrivastava2025} analyse globally convergent shifted QR iteration in finite arithmetic under bounds on eigenvector conditioning and eigenvalue separation, with polylogarithmic precision.

For the eigenpair problem, Armentano, Beltr\'an, B\"urgisser, Cucker, and Shub \cite{ArmentanoEtAl2018} developed a stable polynomial-time algorithm. For Hermitian matrices, Shah \cite{Shah2025} obtains logarithmic precision with nearly matrix multiplication work. The Schur problem considered here requires a basis for the whole space and allows arbitrary complex matrices.

Bai, Demmel, and Gu \cite{BaiDemmelGu1997} developed inverse-free spectral divide and conquer using matrix multiplication and QR factorisation. Schneider's dissertation \cite[Section~1.6.1 and Chapter~6]{Schneider2024} analyses the finite-precision behaviour of several components of inverse-free methods, while Demmel, Dumitriu, and Schneider \cite{DemmelDumitriuSchneider2026} study deflating-subspace algorithms within a common rational framework. The workshop report \cite[Section~3.2]{AmselEtAl2026} gives an account of these results and explains what remains to obtain a complete logarithmic-precision analysis. Our proof of the Schur bound uses Newton iteration with classical inversion, together with an explicit schedule for the requested failure probability.

\section{Arithmetic and conditioning}\label{sec:kernels}
For $M\in\C^{m\times m}$, let $\lambda_1,\ldots,\lambda_m$ be its eigenvalues, counted with multiplicity, and set $\Lambda_M=\diag(\lambda_1,\ldots,\lambda_m)$. Define
\[
 \kappa_V(M)=\inf_{M=V\Lambda_M V^{-1}}\|V\|\|V^{-1}\|,
 \qquad
 \gap(M)=\min_{i\ne j}|\lambda_i-\lambda_j|.
\]
The infimum is over nonsingular diagonalising matrices $V$, with value $+\infty$ if none exists. For a matrix of order one, the gap is $+\infty$. Thus a positive gap implies a simple spectrum. Whenever we have a strict upper bound for $\kappa_V(M)$, we may choose a diagonalising matrix whose condition number lies below that bound.

To obtain logarithmic precision, we need error bounds for multiplication, QR and inversion in which the powers of the condition number are fixed. We use larger dimension factors than necessary to keep the later estimates simple. The kernels and the elementary computations below share a universal constant $C_0\ge2$. We check after the lemma that the specified normalised complex Householder and QR-inverse routines satisfy these bounds with $C_0=2^{20}$.

\begin{lemma}[Classical kernels]\label{lem:kernels}
Let $m\ge1$ be an integer and suppose $C_0m^6u\le1/4$. For matrices with dimensions at most $m$, classical algorithms use $O(m^3)$ operations and have the following properties, with the additional conditioning requirement in the inversion bound.
\begin{enumerate}
\item Matrix multiplication and addition satisfy
\begin{align*}
 \|\fl(XY)-XY\|&\le C_0m^6u\|X\|\|Y\|,\\
 \|\fl(X+Y)-(X+Y)\|&\le C_0m^6u(\|X\|+\|Y\|).
\end{align*}
\item Householder QR of $Y\in\C^{m\times r}$, $0\le r\le m$, including a full unitary completion, returns $\widehat Z,\widehat R$ for which an exactly unitary $Z$ satisfies
\begin{equation}\label{eq:qr}
 \|\widehat Z-Z\|\le C_0m^6u,\qquad
 Y+\Delta=Z\widehat R,\qquad
 \|\Delta\|\le C_0m^6u\|Y\|.
\end{equation}
Here $\widehat Z,Z\in\C^{m\times m}$ and $\widehat R\in\C^{m\times r}$. The matrix $\widehat R$ is upper trapezoidal. Its first $r$ rows, together with the first $r$ columns of $Z$, give the thin factorisation when $Y+\Delta$ has full column rank.
\item Denote by $\INV(Y)$ the inverse computed by Householder QR and triangular solves. For nonsingular square $Y$, it satisfies
\begin{equation}\label{eq:inverse}
 \|\INV(Y)-Y^{-1}\|
 \le C_0m^6u\|Y\|\|Y^{-1}\|^2
\end{equation}
provided $C_0m^6u\kappa_2(Y)\le1/4$.
\end{enumerate}
\end{lemma}
\begin{proof}
Combine the componentwise dot-product bound with $\|X\|_F\le\sqrt m\|X\|$ to obtain the multiplication estimate. Addition is similar. Householder QR gives \eqref{eq:qr}, with a nearby exact unitary product \cite[Chapters~3, 8, and 19]{Higham2002}. We allow a sixth power of the dimension to cover complex arithmetic, the changes of norm and the accumulation of the full unitary factor. For a square input, adjointing the unitary factor matches the interface in \cite[Definition~2.8]{BanksEtAl2023}.

For inversion, we must first know that back substitution is possible. The smallness hypothesis gives $\|Y^{-1}\Delta\|\le C_0m^6u\kappa_2(Y)\le1/4$. Thus $Y+\Delta$ and $\widehat R$ are nonsingular, and back substitution encounters no zero pivot. Let $\widehat X=\INV(Y)$. The backward errors of the columnwise triangular solves, together with the QR error, give
\[
 \|Y\widehat X-I\|\le C_0m^6u(1+\|Y\|\|\widehat X\|).
\]
Multiplication by $Y^{-1}$ and absorption of the term containing $\widehat X$ give $\|\widehat X\|\le2\|Y^{-1}\|$. Substituting this bound and using $\|Y\|\|Y^{-1}\|\ge1$ gives an inverse error bounded by a universal multiple of $m^6u\|Y\|\|Y^{-1}\|^2$. The calculation below gives \eqref{eq:inverse} with the common value $C_0=2^{20}$.
\end{proof}

\paragraph{A concrete kernel constant}
For a complex dot product of length at most $m$, the absolute error is at most $16mu\|x\|\|y\|$. Summing real squares and taking a square root computes a vector norm with relative error at most $16mu$.

Consider first one Householder step. We skip a zero active column. For a nonzero active column $x$, let $v=x+\phi\|x\|e_1$, where $\phi=x_1/|x_1|$ when $x_1\ne0$ and $\phi=1$ otherwise, and let $w=v/\|v\|$. Hats denote the computed counterparts. The phase has error at most $32u$. The unnormalised vector has error at most $64mu\|x\|$, and normalisation gives $\|\widehat w-w\|\le256mu$. The reflector perturbation is therefore at most $2^{11}mu$. Taking Frobenius norms in the dot-product and update error bounds, and retaining exact column zeros, gives a matrix error at most $2^{13}m^2u$ times the current matrix norm.

To bound the QR backward error and the error in the full basis with the same $Z$, take the exact reflector of the current stored active column at each step and compare both recurrences with this sequence. There are at most $m$ reflectors. The smallness condition gives $(1+2^{13}m^2u)^m-1\le2^{14}m^3u$, so $\|\Delta\|\le2^{14}m^3u\|Y\|$ and $\|\widehat Z-Z\|\le2^{14}m^3u$.

For inversion, let $\widehat Z,\widehat R$ be the computed QR factors and $\widehat X$ the computed inverse. Complex division using real and imaginary parts has normwise relative error at most $16u$. Applying the triangular solve bound to each row gives, with entrywise absolute values,
\[
 |\widehat R\widehat X-\widehat Z^*|
 \le64mu\bigl(|\widehat Z^*|+|\widehat R||\widehat X|\bigr).
\]
The QR bounds then imply
\[
 \|Y\widehat X-I\|\le2^{15}m^3u(1+\|Y\|\|\widehat X\|).
\]
Under $2^{20}m^6u\kappa_2(Y)\le1/4$, absorption gives $\|\widehat X\|\le2\|Y^{-1}\|$. The inverse error is at most
$
 3\cdot2^{15}m^3u\|Y\|\|Y^{-1}\|^2.
$
These constants, the multiplication and addition errors, and the elementary shift and Newton errors fit within $C_0m^6u$ with $C_0=2^{20}$.

\section{A stable matrix-sign iteration}\label{sec:sign}
Let $M$ have no eigenvalues on the imaginary axis. Its matrix sign equals $+I$ on the right-half-plane spectral subspace and $-I$ on the left-half-plane spectral subspace. Thus
\begin{equation}\label{eq:projector-intro}
 P_- =\frac{I-\sign(M)}2
\end{equation}
is the spectral projector onto the latter subspace. For a non-normal matrix this projector need not be orthogonal. To obtain a unitary change of basis for the Schur reduction, we use QR to find an orthonormal basis for its range and complete it to a basis for the whole space.

Set $f(z)=(z+z^{-1})/2$ and, for a nonsingular matrix $X$, $f(X)=(X+X^{-1})/2$. The exact Newton iterates are $X_j=f^j(M)$, where $f^j$ denotes $j$-fold composition and $f^0$ is the identity map. We first locate the scalar iterates and bound their distance from the imaginary axis. This will give us one pair of contours for the spectra of all the matrix iterates.

\begin{lemma}[Scalar orbit]\label{lem:scalar}
Let $|z|\le R_0$, $|\Re z|\ge a$, where $R_0\ge1$ and $0<a\le1$, and set
\[
 B=(R_0+1)^2/a\ge4.
\]
For every $j\ge0$,
\begin{equation}\label{eq:scalar-bounds}
 |f^j(z)|\le B,\qquad |\Re f^j(z)|\ge B^{-1}.
\end{equation}
If $\Re z>0$ and $w=(z-1)/(z+1)$, then
\begin{equation}\label{eq:cayley}
 1-|w|\ge2/B,
 \qquad f^j(z)=\frac{1+w^{2^j}}{1-w^{2^j}}.
\end{equation}
\end{lemma}
\begin{proof}
The Cayley transform turns Newton's iteration into squaring. Direct substitution gives the second identity in \eqref{eq:cayley}, while
\[
 1-|w|^2=\frac{4\Re z}{|z+1|^2}\ge\frac4B
\]
gives its first inequality. For $\xi=w^{2^j}$, we have $|\xi|\le|w|<1$, and hence
\[
 \left|\frac{1+\xi}{1-\xi}\right|\le\frac2{1-|w|}\le B,
 \qquad
 \Re\frac{1+\xi}{1-\xi}=\frac{1-|\xi|^2}{|1-\xi|^2}
 \ge\frac{1-|w|}{2}\ge B^{-1}.
\]
The negative half-plane follows from the oddness of $f$ and its iterates.
\end{proof}

The computed iterate can be defective and need not commute with the exact iterate, so we cannot compare them in a common eigenbasis. While the iterates are close, we can enclose both spectra with the same contours. This gives a perturbation bound for all subsequent Newton steps whose constant depends only on the original conditioning and spectral bounds.

\begin{lemma}[Uniform perturbation bound]\label{lem:uniform}
Suppose $M=V\Lambda_M V^{-1}$, with $\Lambda_M$ diagonal and $\kappa_2(V)\le K$, where $K\ge1$. Assume that $|z|\le R_0$ and $|\Re z|\ge a$ for every $z\in\sigma(M)$, with $R_0\ge1$ and $0<a\le1$. Use $B=(R_0+1)^2/a$ from \cref{lem:scalar} and set
\begin{equation}\label{eq:rho-L}
 X_j=f^j(M),\qquad \rho=\frac1{4KB},\qquad L=4096K^2B^6.
\end{equation}
For all nonnegative integers $j,r$ and all matrices $W_1,W_2$ of the same order as $M$ satisfying $\|W_i-X_j\|\le\rho$, the compositions are defined and
\begin{equation}\label{eq:uniform}
 \|f^r(W_1)-f^r(W_2)\|\le L\|W_1-W_2\|.
\end{equation}
\end{lemma}
\begin{proof}
Let $\Gamma$ be the positively oriented boundaries of the two rectangles
\[
 [1/(2B),2B]+\mathrm i[-2B,2B],\qquad
 [-2B,-1/(2B)]+\mathrm i[-2B,2B].
\]
Their combined perimeter is at most $24B$. By \eqref{eq:scalar-bounds}, the spectrum of each $X_j$ lies inside these rectangles, at distance at least $1/(2B)$ from their boundaries. The matrix $V$ diagonalises every $X_j$, even when eigenvalues coalesce during iteration. Therefore
\[
 \|(zI-X_j)^{-1}\|\le2KB\quad (z\in\Gamma).
\]
The definition of $\rho$ and a Neumann series imply
\[
 \|(zI-W_i)^{-1}\|\le4KB\quad (z\in\Gamma).
\]
The same resolvent argument applies along the line segment from $X_j$ to $W_i$. The spectrum avoids $\Gamma$ all along this segment, so no eigenvalue leaves the union of the two rectangles.

Since $\Re f(z)=\Re z(1+|z|^{-2})/2$, each open half-plane is mapped into itself. Every iterate of $f$ is therefore analytic on these half-planes. On $\Gamma$, $|z|\le3B$ and $|\Re z|\ge1/(2B)$. Applying \cref{lem:scalar} with these larger bounds gives, for every $r$,
\[
 |f^r(z)|\le2B(3B+1)^2\le32B^3.
\]
The matrix Cauchy formula and the resolvent identity now yield
$$
 \|f^r(W_1)-f^r(W_2)\|
 \le\frac{24B}{2\pi}\,32B^3(4KB)^2\|W_1-W_2\|\le4096K^2B^6\|W_1-W_2\|.
$$
Since neither the contours nor the bounds on them depend on $r$, the estimate holds for any number of remaining Newton steps.
\end{proof}

\begin{lemma}[Accumulation of local errors]\label{lem:inexact}
Under the hypotheses and notation of \cref{lem:uniform}, let $N\ge1$ be an integer and $e\ge0$. Starting from $Y_0$, form the recurrence below whenever its right-hand side is defined, with error matrices $E_j$ satisfying
\[
 \|Y_0-M\|\le e,\qquad Y_{j+1}=f(Y_j)+E_j,\qquad \|E_j\|\le e\quad(0\le j<N).
\]
If
\begin{equation}\label{eq:e-small}
 e\le\frac{\rho}{4(N+1)L^2},
\end{equation}
then all $N$ iterations are defined and
\begin{equation}\label{eq:additive}
 \|Y_j-X_j\|\le(j+1)Le\qquad(0\le j\le N).
\end{equation}
\end{lemma}
\begin{proof}
The bound holds at $j=0$ because $L\ge1$. Suppose the bound holds through step $s<N$. The iterate $Y_s$ then lies within $\rho/(4L)$ of $X_s$. By \cref{lem:uniform} with $r=1$, $f(Y_s)$ lies within $\rho/4$ of $X_{s+1}$. After the local error is added, $Y_{s+1}$ is still within $\rho/2$. Every remaining composition is therefore defined at both the value before addition and the computed iterate. To isolate the contribution of each local error, we use the identity
\begin{align}\label{eq:telescope}
 Y_j-f^j(M)
 = f^j(Y_0)-f^j(M)+\sum_{s=0}^{j-1}\bigl[f^{j-s-1}(Y_{s+1})
       -f^{j-s-1}(f(Y_s))\bigr].
\end{align}
Each negative term in the sum cancels the preceding positive term by composition. Take $j=s+1$ in \eqref{eq:telescope}. The two arguments in each summand lie within $\rho$ of the corresponding exact iterate, so \cref{lem:uniform} bounds each summand and the initial difference by $Le$. Adding these bounds proves \eqref{eq:additive} at the next step.
\end{proof}

Because the multiplier $L$ is independent of the number of remaining steps, the powers of $K$ and $B$ in \eqref{eq:additive} stay fixed as $N$ increases.

\begin{corollary}[Sign computation]\label{cor:sign}
Suppose $M\in\C^{m\times m}$ satisfies the hypotheses of \cref{lem:uniform}. Then $\sign(M)$ can be approximated with forward error at most $\varepsilon\in(0,1)$ using
$
 O\bigl(m^3\log(mKB/\varepsilon)\bigr)
$
operations and $O(\log(mKB/\varepsilon))$ mantissa bits. Upper bounds for $m,K,B$ suffice as inputs to the choice of precision and iteration count.
\end{corollary}
\begin{proof}
Take
\begin{equation}\label{eq:N}
 j_*=\left\lceil\log_2\frac{8K}{\varepsilon}\right\rceil,\qquad
 N=\left\lceil\log_2(Bj_*)\right\rceil.
\end{equation}
By \eqref{eq:cayley}, $|w|^{2^N}\le \exp(-2j_*)\le2^{-j_*-1}\le\varepsilon/(16K)$. The scalar distance of $f^N(z)$ from the appropriate sign is at most $2|w|^{2^N}/(1-|w|^{2^N})$. Multiplying by $K$ bounds $\|X_N-\sign(M)\|$ by $\varepsilon/4$. Because each step squares the Cayley variable, the fixed count grows only doubly logarithmically with $K/\varepsilon$ for a fixed spectral enclosure.

For $\|Y_j-X_j\|\le\rho$, \eqref{eq:scalar-bounds} gives $\|X_j\|,\|X_j^{-1}\|\le KB$. Another Neumann series then gives
$
 \|Y_j\|,\|Y_j^{-1}\|\le2KB.
$
Implement the update with the classical inverse in \cref{lem:kernels}. Its local error, together with addition and scaling errors, is at most
\begin{equation}\label{eq:local-e}
 e=C_0m^6u(KB)^3
\end{equation}
with the same $C_0=2^{20}$ as in the kernel bounds above. After halving, the inverse error is at most $3\cdot2^{17}m^3u(KB)^3$, and the halved addition error is at most $3\sqrt m\,uKB$. Their sum is below $C_0m^6u(KB)^3$ for $KB\ge4$. Binary halving is exact. Choose $u$ so that this bound is at most both quantities
\[
 \frac{\rho}{4(N+1)L^2},\qquad
 \frac{\varepsilon}{2(N+1)L}.
\]
The inverse smallness condition in \cref{lem:kernels} must also hold. Each upper bound for the roundoff is an inverse polynomial in $m,K,B,\varepsilon^{-1},N$ with fixed exponents. Since $N=O(\log(KB/\varepsilon))$, $O(\log(mKB/\varepsilon))$ mantissa bits suffice for all of them. The initial error includes the rounding of the input matrix and any scalar shift. By \cref{lem:inexact}, the iterates stay in the neighbourhood where the inverse and local error estimates are valid, and the accumulated error is at most the remaining $\varepsilon/2$.
\end{proof}

\section{Conditioning of perturbed matrices and diagonal blocks}\label{sec:conditioning}
Subdivision perturbs the matrix through rounding and the removal of lower-left blocks. To use sign iteration on the next diagonal blocks, we must bound their eigenvector conditioning after these perturbations. We first bound the conditioning of the perturbed full matrix. We then pass to its diagonal blocks by compressing its spectral projectors.

\begin{lemma}[A robust eigenvector basis]\label{lem:robust}
Let $M=V\Lambda_M V^{-1}\in\C^{n\times n}$, with $\Lambda_M$ diagonal. Assume that $M$ has simple spectrum, $\kappa_2(V)\le K$ and $\gap(M)\ge g$, where $K\ge1$ and $g>0$. If $F\in\C^{n\times n}$ satisfies
$
 \|F\|\le\frac{g}{16\sqrt n K},
$
then $M+F$ has simple spectrum. Its eigenvalues match uniquely to those of $M$, each within $K\|F\|$, and
\begin{equation}\label{eq:robust}
 \kappa_V(M+F)\le3K.
\end{equation}
\end{lemma}
\begin{proof}
In the original eigenbasis the perturbation is $E=V^{-1}FV$, with $\|E\|\le K\|F\|$. The Bauer--Fike inclusion follows directly from the resolvent of $\Lambda_M$: every eigenvalue of $\Lambda_M+tE$ lies within $t\|E\|$ of an original eigenvalue. The discs of radius $\|E\|$ are pairwise disjoint. Along the homotopy $t\in[0,1]$, each therefore retains exactly one eigenvalue, counted with multiplicity.

For the eigenvalue $\mu_i$ near $\lambda_i$, delete the $i$th row and column of $\Lambda_M-\mu_i I$. The remaining diagonal entries have modulus at least $g-\|E\|$, and the perturbed complementary system has inverse norm at most $2/g$. If an eigenvector had zero $i$th coordinate, this invertible system would force all its coordinates to vanish. We can therefore normalise it as $e_i+y_i$, where $(y_i)_i=0$. The complementary equations give
$
 \|y_i\|\le2\|E\|/g.
$
If $Y$ has columns $y_i$, then $\|Y\|\le\|Y\|_F\le2\sqrt n\|E\|/g\le1/8$. Thus $V(I+Y)$ is an eigenvector basis and its condition number is at most $K(1+1/8)/(1-1/8)<3K$.
\end{proof}

\begin{lemma}[Diagonal-block conditioning]\label{lem:block}
Let $S$ be upper block triangular, with simple spectrum and $\kappa_V(S)\le K$. Every diagonal block $C$ of order $m$ satisfies
\begin{equation}\label{eq:block-condition}
 \kappa_V(C)\le mK.
\end{equation}
\end{lemma}
\begin{proof}
For a simple eigenvalue $\lambda_i$ of $S$, its spectral projector is $p_i(S)$, where
\[
 p_i(z)=\prod_{j\ne i}\frac{z-\lambda_j}{\lambda_i-\lambda_j}.
\]
Its norm is at most $K$. Applying a polynomial to an upper block triangular matrix applies the same polynomial to each diagonal block. Thus, if $\lambda_i\in\sigma(C)$, its projector $p_i(C)$ is the diagonal block of $p_i(S)$ corresponding to $C$ and has norm at most $K$.

The simple spectrum of $S$ implies that of $C$. Normalise the right eigenvectors of $C$ to unit norm and form the matrix $W$ with these columns. Then $\|W\|\le\sqrt m$. Each rank-one spectral projector is the product of a column of $W$ and the corresponding row of $W^{-1}$. Since that column has unit norm, the projector norm equals the row norm. Hence $\|W^{-1}\|\le\sqrt m K$, proving \eqref{eq:block-condition}. If the infimum defining $\kappa_V$ is not attained, take the limit along a sequence of diagonalising bases.
\end{proof}

The polynomial-projector argument works even for an interior block, which need not occupy an invariant coordinate subspace of the full matrix. Applying the lemma to the full matrix at each level introduces the factor $m$ only once for each current block.

\section{Parameter selection and finite randomness}\label{sec:random}
For $n=1$, return $Q=1$ and the input scalar rounded with $b\ge\max\{8,\lceil\log_2(1/\delta)\rceil+2\}$. The resulting error is deterministically below $\delta$. For the rest of the construction, assume $n\ge2$.

Round the requested $\delta,p$ down to powers of two in $(\delta/2,\delta]$ and $(p/2,p]$, respectively. In this section and the proof below, $\delta,p$ denote the rounded values. This changes $\log(n/(\delta p))$ by at most a universal additive constant. We first derive sufficient conditions for these requests, then choose a schedule that fixes the sampling, grid and arithmetic parameters before any draws. The same conditions also allow less conservative choices at fixed confidence. In the bounds below, $\sigma$ is a rational upper bound for $\sqrt n$ with $\sqrt n\le\sigma\le2\sqrt n$. The implementation uses 64 fractional binary places for this bound.

\subsection{A parameterised finite sampler}
The sampler approximates one standard complex Gaussian at a time. Choose an integer $\ell\ge19$ and set $a=2^{-\ell}$, $W=\lceil\sqrt\ell\rceil$ and $V_{\rm G}=nW$. The cutoff $a$ controls how often we reject a draw, and $W$ bounds the reference Gaussian on the retained event. For a requested entrywise accuracy $0<\zeta\le1/4$, set
\begin{equation}\label{eq:sampler-parameters}
 m_U=\left\lceil\log_2\frac{64}{\zeta a^2}\right\rceil,
 \qquad J=32+\left\lceil\log_2\frac{2^{16}(\ell+1)}{\zeta a}\right\rceil.
\end{equation}
Choose $\widehat U_1,\widehat U_2$ independently and uniformly on the mesh of width $2^{-m_U}$, using $m_U$ bits per endpoint. If the radial endpoint $\widehat U_1$ lies outside $[2a,1-2a]$, reject the draw without retrying. For the coupling, take $U_i$ uniformly in the cell with left endpoint $\widehat U_i$, independently for the two coordinates. The reference Box--Muller variable is
\[
 G=\sqrt{-\log U_1}\bigl(\cos(2\pi U_2)+\mathrm i\sin(2\pi U_2)\bigr),
\]
with density $\pi^{-1}e^{-|z|^2}$ \cite{BoxMuller1958}. To compute $\widehat G$, evaluate the formula at $\widehat U_1,\widehat U_2$ with truncated series and rounded arithmetic. Reduce the logarithm using $\log(2^kv)=k\log2+\log v$, with $k\in\mathbb Z$ and $1\le v\le2$, and use
\[
 \log v=2\sum_{j\ge0}\frac{w^{2j+1}}{2j+1},
 \qquad w=\frac{v-1}{v+1},\quad |w|\le1/3.
\]
Use the same series for $\log2$, the Machin identity $\pi=16\arctan(1/5)-4\arctan(1/239)$, and Taylor series for sine and cosine. Each series has $J$ terms, generated by recurrence.

\begin{lemma}[Finite sampler]\label{lem:sampler}
Suppose the integer coefficients and uniforms are exactly representable and
\begin{equation}\label{eq:sampler-u}
 u\le\frac{\zeta a}{2^{60}(\ell+1)^2(J+1)^2}.
\end{equation}
For at most $n^3$ complex entries, the sampler can be coupled to independent standard complex Gaussians such that, outside an event of probability at most $6n^3a$, no rejection occurs and
\begin{equation}\label{eq:sample-bound}
 |\widehat G-G|\le\zeta,\qquad |G|\le W.
\end{equation}
Its work per entry is $O(J+m_U)$, including random-bit generation. Every loop is deterministically capped.
\end{lemma}
\begin{proof}
Couple a mesh endpoint to a continuous uniform in its cell, of width $\Delta=2^{-m_U}\le\zeta a^2/64$. The set leading to endpoint rejection has probability at most $6a$. On its complement the radial endpoint and its continuous reference lie in $[a,1-a]$. The radial and angular derivatives of the complex Box--Muller map are bounded by $a^{-2}$ there, so mesh error is at most $2\Delta/a^2\le\zeta/32$. Also $|G|\le\sqrt{\ell\log2}<W$.

The logarithm and arctangent tails are bounded by $2^{-J}$ before range reduction. For arguments of modulus at most seven, the sine and cosine tails, starting at $J\ge32$, are smaller than $2^{-J}$. For example, use $k!\ge(k/3)^k$ and bound the remaining tail geometrically. Equation~\eqref{eq:sampler-parameters} gives $2^{-J}\le\zeta a/(2^{48}(\ell+1))$. Range reduction, the square root, and the final products therefore give total truncation error at most $2^8(\ell+1)2^{-J}/a\le2^{-40}\zeta$.

We can bound scalar roundoff without a matrix-kernel assumption. The reduced logarithm has error at most $16(J+1)u$, range reduction at most $32(\ell+1)(J+1)u$, and the Machin evaluation at most $2^{12}(J+1)^2u$. For each trigonometric series, the sum of the absolute values of its Taylor terms is at most $e^7<2^{11}$. The product recurrence and summation bounds give error at most $2^{20}(J+1)^2u$ per function. Since the accepted radial endpoint $\widehat U_1$ satisfies $-\log\widehat U_1\ge2a$, the combined log error is less than $a/4$. The square root is therefore positive and amplifies absolute error by at most $a^{-1}$. Propagating these bounds through the angle and products gives at most $2^{24}(\ell+1)^2(J+1)^2u/a\le2^{-36}\zeta$. These estimates use $(J+1)u\le2^{-60}$, which follows from \eqref{eq:sampler-u}. Adding mesh, truncation and rounding errors proves \eqref{eq:sample-bound}. A conditional union bound over the capped number of draws proves the failure bound, including adaptively requested range samples.
\end{proof}

\subsection{Regularisation with a requested confidence}
Let $G\in\C^{n\times n}$ have independent standard complex Gaussian entries and let $\widehat G$ be its coupled sampled matrix. Set $G_n=G/\sqrt n$. Choose a perturbation scale $0<\gamma<1/2$ and a norm cutoff $M_{\rm cut}>0$ with $\gamma M_{\rm cut}\le\delta/4$. We seek an eigenvector condition number below $t$ and a gap above $r$, where $t\ge1$ and $0<r\le1$. The regularisation estimate bounds the probability that the Gaussian norm meets the cutoff but either spectral target fails.
\begin{equation}\label{eq:BGVKS}
 \Prob\bigl\{(\kappa_V(A+\gamma G_n)\ge t\text{ or }\gap(A+\gamma G_n)\le r)
                    \text{ and }\|G_n\|\le M_{\rm cut}\bigr\}\le \frac{9n^3}{\gamma^2t^2}
       +\frac{576t^8r^6n^8}{\gamma^8}.
\end{equation}
This estimate follows from the proof of \cite[Theorem~3.6]{BanksEtAl2023}, with the norm cutoff used here. On $\|G_n\|\le M_{\rm cut}$, $\|A+\gamma G_n\|\le1+\delta/4<2$, so the same radius-three spectral disc and its $r/2$-net apply. The eigenvector-condition estimate on that disc gives $9n^3/(\gamma^2t^2)$. The net has at most $144/r^2$ points, and the two-eigenvalue estimate at each point is $4(trn/\gamma)^8$. Their product gives the second term.

Let $\omega$ bound the probability that $\|G_n\|>M_{\rm cut}$ while all initial sampler uniforms are retained. For the choice $M_{\rm cut}=V_{\rm G}$, \cref{lem:sampler} gives $\omega=0$: on the retained event, $\|G_n\|\le\|G_n\|_F\le\sqrt nW\le V_{\rm G}$. 

To pass from the Gaussian perturbation to the computed matrix, allow a formation error $\nu>0$. We require
\begin{equation}\label{eq:regularisation-budgets}
 \nu\le\min\left\{\frac\delta{48},\frac{r}{32\sigma t}\right\},
 \qquad K_0=3t,\qquad g_0=r/2.
\end{equation}
The sampling and roundoff bounds verified in \cref{sec:completion} ensure
$\|M_0-(A+\gamma G_n)\|\le\nu$ for the computed regularisation
$M_0=\fl(A+\gamma\widehat G/\sqrt n)$, including input rounding.
When $\|G_n\|\le M_{\rm cut}$, $\kappa_V(A+\gamma G_n)<t$ and $\gap(A+\gamma G_n)>r$, \cref{lem:robust} then yields
\begin{equation}\label{eq:M0}
 \|M_0-A\|\le\delta/3,\quad \|M_0\|\le2,
 \quad\kappa_V(M_0)\le K_0,\quad\gap(M_0)\ge g_0.
\end{equation}
The exceptional probability is at most the two terms in \eqref{eq:BGVKS}, plus $\omega$ and the sampler failure already counted once.

\subsection{Grid, local budgets and precision}
We bisect each coordinate at most $H$ times. Thus $D$ is the cap on the total number of levels, and $h$ is the spacing of the finest grid. We also choose the bound $K$ for the eigenvector conditioning of a current block, justified in \cref{sec:global}. Set
\begin{equation}\label{eq:grid-parameters}
 H=\lceil\log_2(32/g_0)\rceil,\quad D=2H,\quad h=8\,2^{-H},\quad K=3nK_0.
\end{equation}
Choose a grid margin $0<\eta\le h/(16n)$ and let $M_G$ be the number of possible offsets per coordinate, a positive power of two. To control the perturbations from subdivision, set the total allowance $\mu$ and the allowance $d$ per level by
\begin{equation}\label{eq:level-budgets}
 \mu=\min\left\{\delta/12,g_0/(32\sigma K_0),\eta/(4K_0)\right\},
 \qquad d=\min\{\mu/D,\eta/(16nK)\}.
\end{equation}
Choose offsets $\alpha_x,\alpha_y$ independently and uniformly from $\{jh/M_G:j=0,\ldots,M_G-1\}$. They translate the root square to $[\alpha_x-4,\alpha_x+4]\times[\alpha_y-4,\alpha_y+4]$. We bisect its coordinates alternately. All vertical and horizontal cut coordinates then belong to $\alpha_x+h\mathbb Z$ and $\alpha_y+h\mathbb Z$, respectively.

\begin{lemma}[Simultaneous grid margin]\label{lem:grid}
Conditional on the initial matrix satisfying \eqref{eq:M0}, the probability that some initial eigenvalue is within $\eta$ of a vertical or horizontal grid line is at most
\begin{equation}\label{eq:grid-failure}
 \pi_G=4n\eta/h+4n/M_G.
\end{equation}
\end{lemma}
\begin{proof}
Fix one eigenvalue and one coordinate. Modulo the grid spacing, the offsets that bring a line within $\eta$ of this coordinate form an arc of length $2\eta$ on a circle of circumference $h$. The fraction of mesh points in the arc exceeds its fraction of the circumference by at most $2/M_G$. Summing over both coordinates and all $n$ eigenvalues gives the bound. Choosing the offsets requires no eigenvalue computation.
\end{proof}

At a nontrivial split, we seek a lower bound $s\in(0,1)$ for the smallest singular value of the product of the exact spectral projector with a fresh Gaussian matrix. In \cref{sec:subspace}, $\varepsilon_Y$ bounds the perturbation of this product, $\tau$ is the accuracy of the computed projector, and $\zeta$ is the common entrywise sampling accuracy. Set
\begin{align}\label{eq:local-budgets}
 \varepsilon_Y&=sd/64,\qquad
 \tau=\min\{1/(16n),\varepsilon_Y/(4V_{\rm G})\},\\
 \zeta&=\min\{1/4,\varepsilon_Y/(16Kn),\nu/(8\gamma\sigma)\}.\nonumber
\end{align}
Use $\zeta$ in \eqref{eq:sampler-parameters}. For the sign computation, use the enclosure and Newton count from \cref{lem:uniform,cor:sign} with the following bounds.
\begin{align}\label{eq:audit-sign}
 B&=162/\eta,\quad L=4096K^2B^6,\quad\rho=(4KB)^{-1},\\
 j_*&=\lceil\log_2(8K/\tau)\rceil,\quad N=\lceil\log_2(Bj_*)\rceil.\nonumber
\end{align}
The count $N$ is fixed before any iteration. We express the matrix kernel error bounds in terms of $q=C_0n^6u$ and require every inequality in \cref{tab:budgets}.

\begin{table}[tbp]
\caption{Sufficient bounds for the common arithmetic precision, with $q=C_0n^6u$.}
\label{tab:budgets}
\centering\small
\begin{tabular}{ll}
\toprule
Purpose & Upper bound for $q$\\
\midrule
Kernel smallness & $1/4$\\
Initial formation & $\nu/[64\sigma(1+\gamma W)]$\\
Range product and QR & $\varepsilon_Y/(32K V_{\rm G})$\\
Masked level update & $d/(128\sigma)$\\
Basis accumulation & $\delta/(1024D)$\\
Trace rank & $1/(64nK)$\\
Projector formation & $\tau/[4(K+1)]$\\
Inverse smallness & $1/[16(KB)^2]$\\
Sign stability & $1/[2^{31}(N+1)K^8B^{16}]$\\
Sign accuracy & $\tau/[2^{16}(N+1)K^5B^9]$\\
\bottomrule
\end{tabular}
\end{table}

In addition impose \eqref{eq:sampler-u} and
\begin{equation}\label{eq:representation-bits}
 b\ge m_U+2,\qquad
 b\ge2\lceil\log_2(2J+2)\rceil+2,\qquad
 b\ge H+\log_2M_G+4.
\end{equation}
These three bounds ensure that we can represent the uniform mesh samples, integer series coefficients and dyadic grid coordinates exactly. If $q_*$ is the minimum of the table's right-hand sides, an explicit sufficient precision is
\begin{equation}\label{eq:precision-choice}
 b_{\rm dir}=\max\left\{
 \left\lceil\log_2\frac{C_0n^6}{q_*}\right\rceil,
 \left\lceil\log_2\frac{2^{60}(\ell+1)^2(J+1)^2}{\zeta a}\right\rceil,
 \text{the bounds in \eqref{eq:representation-bits}}\right\}.
\end{equation}
The implementation's safe mode uses $b_{\rm dir}+32$, while direct mode uses $b_{\rm dir}$. Neither mode assumes an affine bound.

\subsection{An explicit confidence-dependent choice}
For the rounded confidence request, let $\beta=p/8$. Define
\begin{align}\label{eq:explicit-confidence}
 \ell&=\max\{19,\lceil\log_2(8n^3/\beta)\rceil\},\quad
 \gamma=\delta\,2^{-\lceil\log_2(16V_{\rm G})\rceil},\\
 t&=2^{\lceil\frac12\log_2(9n^3/(\gamma^2\beta))\rceil},\quad
 r=2^{\lfloor\frac16\log_2(\beta\gamma^8/(576t^8n^8))\rfloor},\nonumber\\
 \eta&=h\beta/(8n),\quad
 M_G=2^{\lceil\log_2(8n/\beta)\rceil},\quad
 s=2^{\lfloor\frac12\log_2(\beta/n^3)\rfloor}.\nonumber
\end{align}
Take $\nu$ equal to the upper bound in \eqref{eq:regularisation-budgets}. The order of the choices matters. First $\ell$ determines $a,W,V_{\rm G}$, and then we choose $\gamma,t,r$. These give the regularisation and grid parameters, after which \eqref{eq:local-budgets} determines the local accuracies. We can then choose the sampler, sign and precision parameters. Integer comparisons give conservative dyadic bounds for the roots and logarithms without evaluating an approximate logarithm. With these choices, $\gamma V_{\rm G}\le\delta/16$, both terms of \eqref{eq:BGVKS} are at most $\beta$, $\pi_G\le\beta$, $n^3s^2\le\beta$, and $6n^3a<\beta$.

\section{From a projector to an invariant subspace}\label{sec:subspace}
Let $C\in\C^{m\times m}$ be the exact matrix represented by a current stored block, with $m\le n$. For a vertical cut at $c\in\R$, the exact shifted matrix is $M=C-cI$. For a horizontal cut it is $M=-\mathrm iC-cI$. In either case the computation starts from the rounded shift. We need three bounds on $C$, which follow from the invariant in \cref{sec:global}.
\begin{equation}\label{eq:block-assumptions}
 \|C\|\le3,\qquad\kappa_V(C)\le K,\qquad
 \operatorname{dist}(\sigma(C),\text{cut line})\ge\eta/2.
\end{equation}
The grid gives $|c|\le4+h$, so the shifted matrix has norm at most eight and \cref{lem:scalar} applies with $B=162/\eta$. The exact sign projector $P$ selects the eigenvalues to the left of a vertical cut or below a horizontal cut. Let $\widehat{\sign}(M)$ denote the result of exactly $N$ floating-point Newton steps from \eqref{eq:audit-sign}, and form $\widetilde P=\fl((I-\widehat{\sign}(M))/2)$. The precision bounds verified in \cref{sec:completion} give
\begin{equation}\label{eq:projector-accuracy}
 \|\widetilde P-P\|\le\tau,\qquad P=(I-\sign(M))/2.
\end{equation}
The trace of $P$ is its rank. We sum the real parts of the diagonal entries of $\widetilde P$ and take the nearest integer $r_P$. The projector error changes the sum by at most $n\tau\le1/16$. Summation and conversion contribute at most $4nKq\le1/16$. Thus the rank is correct and the distance test at $1/3$ cannot reject. A rank of zero or $m$ places the whole spectrum in one half-cell, so only the cell assignment changes.

\begin{lemma}[Range extraction for an oblique projector]\label{lem:range}
If $P^2=P\in\C^{m\times m}$ has rank $1\le r_P\le m$, and $\Omega$ is an $m\times r_P$ matrix with independent standard complex Gaussian entries, independent of $P$, then, for every $\theta>0$,
\begin{equation}\label{eq:smallball}
 \Prob\{\sigma_{r_P}(P\Omega)<\theta\}\le r_P^2\theta^2.
\end{equation}
\end{lemma}
\begin{proof}
In an orthonormal basis adapted to its range, $P$ has the form
$$
 \begin{pmatrix}I&F\\0&0\end{pmatrix}.
$$
Thus all its nonzero singular values are at least one. A thin SVD $P=U_{r_P}\Sigma V_{r_P}^*$ gives a square standard complex Gaussian $G=V_{r_P}^*\Omega$, so $\sigma_{r_P}(P\Omega)\ge\sigma_{r_P}(G)$.

If $\sigma_{r_P}(G)<\theta$, choose a unit vector $v$ with $\|Gv\|<\theta$ and an index $j$ with $|v_j|\ge1/\sqrt{r_P}$. The $j$th column of $G$ is within $\sqrt{r_P}\,\theta$ of the span of its other columns. Conditional on those other columns, which span an $(r_P-1)$-dimensional space almost surely, the orthogonal component of that column is a standard complex Gaussian scalar. Its probability of having modulus below $\sqrt{r_P}\,\theta$ is $1-e^{-r_P\theta^2}\le r_P\theta^2$. Taking the union over $j$ proves the result.
\end{proof}

At a nontrivial split, condition on the history up to the sample request. The sampler returns $\widehat\Omega\in\C^{m\times r_P}$ coupled to a fresh standard complex Gaussian matrix $\Omega$ independent of that history. We compare the computed product $\fl(\widetilde P\widehat\Omega)$ and its full Householder QR completion $\widehat Z$ with the exact product $Y=P\Omega$. On the sampler event, $\|\Omega\|\le V_{\rm G}$. We also require the range event
\begin{equation}\label{eq:range-good}
 \sigma_{r_P}(P\Omega)\ge s.
\end{equation}
The projector error contributes at most $\tau V_{\rm G}\le\varepsilon_Y/4$, the sample error at most $Kn\zeta\le\varepsilon_Y/16$, and the product and QR bounds in \cref{tab:budgets} absorb their arithmetic errors and cross terms. Thus the exact unitary comparison $Z$ in \cref{lem:kernels} factors $Y+E$, where
\begin{equation}\label{eq:range-error}
 Y=P\Omega,\qquad Y+E=Z\widehat R,\qquad\|E\|\le\varepsilon_Y.
\end{equation}
Since $\|E\|\le sd/64<s/2$, the matrix $Y+E$ has full column rank. The two ranges now have the same dimension, and the largest sine of their principal angles satisfies
\begin{equation}\label{eq:angle}
 \|\sin\Theta(\range Y,\range(Y+E))\|
 \le\frac{\|E\|}{\sigma_{r_P}(Y)-\|E\|}\le d/32.
\end{equation}
Indeed, a unit vector in $\range(Y+E)$ has coefficients of norm at most $1/\sigma_{r_P}(Y+E)$. Orthogonal projection onto the complement of $\range Y$ removes $Y$ and leaves only $E$. Rotating the full frame in its principal planes produces an exactly unitary $Z_0$, with its first $r_P$ columns spanning $\range P$, such that
\begin{equation}\label{eq:direct-rotation}
 \|Z-Z_0\|\le3d/32.
\end{equation}
Partition after the first $r_P$ coordinates. Since $\range P$ is $C$-invariant, $(Z_0^*CZ_0)_{21}=0$. Together with $\|C\|\le3$, this gives
\begin{equation}\label{eq:discard}
 \|(Z^*CZ)_{21}\|\le 2\|C\|\|Z-Z_0\|\le9d/16.
\end{equation}
By \cref{lem:range}, each range failure has conditional probability at most $r_P^2s^2$. We apply this bound before conditioning on sampler acceptance.

\section{Spectral subdivision by levels}\label{sec:algorithm}
The stored state consists of an upper block triangular matrix $S$, a partition into contiguous diagonal blocks with assigned grid cells, and an accumulated basis $\widehat Q$. In \cref{alg:schur}, the local splits at each level give a single change of basis for the full matrix.

\begin{algorithm}[t]
\caption{Schur decomposition with the explicit sufficient schedule.}
\label{alg:schur}
\small
\begin{tabbing}
\hspace{1.2em}\=\hspace{1.2em}\=\hspace{1.2em}\=\kill
\textbf{procedure} \textsc{Schur}($A\in\C^{n\times n},\delta,p$), returning $(Q,T)$\\
\>\textbf{if} $n=1$ \textbf{then return} $(1,\fl(A))$ at the scalar precision in \cref{sec:random}.\\
\>Round $\delta,p$ down and select the schedule \eqref{eq:explicit-confidence} as in \cref{sec:random}.\\
\>Set $b=b_{\rm dir}$ (or $b_{\rm dir}+32$ in safe mode) from \eqref{eq:precision-choice}, and $R$ from \eqref{eq:R}.\\
\>Draw $\widehat G\in\C^{n\times n}$ using the finite sampler of \cref{lem:sampler}.\\
\>$S\leftarrow M_0=\fl(A+\gamma\widehat G/\sqrt n)$, $\widehat Q\leftarrow I$.\\
\>Draw the two independent grid offsets and assign one block to the root square.\\
\>\textbf{for} $k=1,\ldots,D$ \textbf{do}\\
\>\>\textbf{for each} current nonscalar diagonal block $C$ of $S$ \textbf{do}\\
\>\>\>Call \textsc{SplitBlock}($C,\text{cell},k$) and keep its rank, basis and cells.\\
\>\>\textbf{end for}\\
\>\>Assemble $\widehat Y_k$ from these bases and identities on scalar blocks.\\
\>\>Apply the full matrix and basis updates \eqref{eq:level-update}.\\
\>\>Adopt the returned cells and split each parent of order $m$ after coordinate $r_P$\\
\>\>\quad when $0<r_P<m$. Retain scalar blocks and their cells.\\
\>\>Set the new lower-left blocks and all earlier structural zeros exactly to zero.\\
\>\>\textbf{if} all blocks are scalar \textbf{then return} $(\widehat Q,S)$.\\
\>\textbf{end for}\\
\>\textbf{return} $(I,0)$.\\[5pt]
\textbf{procedure} \textsc{SplitBlock}($C,\text{cell},k$), with $C\in\C^{m\times m}$\\
\>Choose a vertical cut for odd $k$ and a horizontal cut for even $k$.\\
\>Let $c$ be the exact midpoint of that cell coordinate and form the two half-cells.\\
\>$Y\leftarrow\fl(C-cI)$ for a vertical cut, or $Y\leftarrow\fl(-\mathrm iC-cI)$ otherwise.\\
\>\textbf{for} $j=1,\ldots,N$ \textbf{do}\\
\>\>$Y\leftarrow\fl((Y+\INV(Y))/2)$.\\
\>\textbf{end for}\\
\>$\widetilde P\leftarrow\fl((I-Y)/2)$, $v\leftarrow\fl(\sum_{i=1}^m\Re\widetilde P_{ii})$.\\
\>Let $r_P$ be the nearest integer to $v$.\\
\>Abort if $r_P\notin\{0,\ldots,m\}$ or the computed trace distance exceeds $1/3$.\\
\>\textbf{if} $r_P=0$ \textbf{then return} $(r_P,I_m,\text{right/upper half-cell})$.\\
\>\textbf{if} $r_P=m$ \textbf{then return} $(r_P,I_m,\text{left/lower half-cell})$.\\
\>Draw a fresh $\widehat\Omega\in\C^{m\times r_P}$ using the finite sampler of \cref{lem:sampler}.\\
\>Compute the full Householder QR basis $\widehat Z$ of $\fl(\widetilde P\widehat\Omega)$.\\
\>\textbf{return} $(r_P,\widehat Z,\text{left/lower then right/upper half-cells})$.
\end{tabbing}
Apply all prescribed sampler and kernel guards, including the global $n^3$ Gaussian-entry cap. Any abort in either routine terminates \textsc{Schur} with $(Q,T)=(I,0)$, without retrying.
\end{algorithm}

Let $\widehat Y_k$ be the block diagonal matrix of computed bases and identity blocks. Applying this change of basis to $S$ and accumulating it in $\widehat Q$ gives the two updates whose errors we analyse.
\begin{equation}\label{eq:level-update}
 S\leftarrow\fl(\widehat Y_k^*S\widehat Y_k),\qquad
 \widehat Q\leftarrow\fl(\widehat Q\widehat Y_k).
\end{equation}
The error bounds hold at every completed level, so we can stop as soon as all blocks are scalar.

\subsection{Termination and exceptional arithmetic}
The sampler, Newton iteration and matrix kernels have the fixed loop lengths selected above. A subdivision tree has at most $n-1$ nontrivial splits. Each range sample uses at most $n(n-1)$ entries, so, including regularisation, the total is at most $n^2+n(n-1)^2\le n^3$. Thus the Gaussian-entry cap cannot reject a valid subdivision tree.

Use the zero-column Householder branch. Abort on a zero divisor, a nonpositive logarithm argument or a nonpositive computed value of $-\log\widehat U_1$, a rejected radial endpoint, an invalid nearest-integer rank, or a computed trace distance exceeding $1/3$. To bound large quantities, choose a power of two
\begin{equation}\label{eq:R}
 R=2^{\left\lceil\log_2\max\{2^{40},C_0^{100},n,K,B,V_{\rm G},
 \nu^{-1},\eta^{-1},s^{-1},\zeta^{-1},\tau^{-1},d^{-1},D,J\}\right\rceil}.
\end{equation}
After each matrix kernel, reject if any component has magnitude at least $R^{100}/4$. The arithmetic model excludes hardware exponent overflow and underflow. To check that the magnitude test cannot reject on the good event, follow the induction in \cref{prop:invariant} with the test suppressed. Newton iterates, inverse outputs, projectors and triangular factors from inverse calls have norms at most $4KB$. Range products have norm at most $25KV_{\rm G}/8$, and their triangular QR factors at most $125KV_{\rm G}/32<4KV_{\rm G}$. The remaining checked outputs have norm below ten. The formed and shifted inputs are bounded by two and ten, respectively. Full transformation products are bounded by $2(1+q)^4<5$, basis products by four, and full QR bases by $1+q$. Every checked entry is therefore bounded by
$
 10+4KB+4KV_{\rm G}\le18R^2<R^{100}/4.
$
It follows that every checked matrix passes the test along this trajectory. We need no magnitude bound for the unchecked scalar accumulators. The exceptional return in \cref{alg:schur} is exactly triangular and keeps the worst-case operation bound deterministic.

\section{Global backward error and spectral labels}\label{sec:global}
Let $S_k,\widehat Q_k$ be the stored matrices after $k$ completed levels, with $S_0=M_0$ and $\widehat Q_0=I$. At level $k\ge1$, let $Z_k$ be the exact unitary block diagonal comparison for the full change of basis $\widehat Y_k$. Define the accumulated exact basis by $W_0=I$ and $W_k=Z_1\cdots Z_k$. The algorithm uses the stored bases, and we use $Z_k,W_k$ to track the errors.

\begin{proposition}[Admissible schedules]\label{prop:invariant}
Suppose \eqref{eq:M0}, the simultaneous grid margin and all retained-sampler and range events hold. Under \eqref{eq:regularisation-budgets}--\eqref{eq:precision-choice}, at every completed level $0\le k\le D$ there is a perturbation $F_k$ such that
\begin{equation}\label{eq:invariant}
 S_k=W_k^*(M_0+F_k)W_k,\qquad\|F_k\|\le k d\le\mu.
\end{equation}
Each eigenvalue of each current block matches uniquely to an initial eigenvalue of $M_0$ in its assigned cell. Every block is scalar by level $D$. At termination both inequalities in \eqref{eq:target} hold.
\end{proposition}
\begin{proof}
The initial state satisfies the invariant. The construction for \eqref{eq:M0} gives an eigenvector basis for $M_0$ with condition number at most $9t/7<K_0$. Under \eqref{eq:M0} alone, the infimum is attained by compactness, using a norm-one minimising sequence with bounded inverses. If the invariant holds after $k$ levels, then $\|F_k\|\le\mu\le g_0/(32\sigma K_0)$, so \cref{lem:robust} gives
\begin{equation}\label{eq:global-condition}
 \kappa_V(M_0+F_k)\le3K_0,\qquad
 |\lambda_i(M_0+F_k)-\lambda_i(M_0)|\le K_0\mu\le\eta/4.
\end{equation}
The current gap is at least $g_0/2$. Also $\|S_k\|\le\|A\|+\delta/3+\mu<2$. By \cref{lem:block}, each diagonal block has condition number at most $3nK_0=K$. The eigenvalue displacement reduces the initial grid margin by at most $\eta/4$, proving \eqref{eq:block-assumptions} for every block about to be split.

The discarded blocks from distinct parents form a direct sum, so their combined norm is at most $9d/16$ by \eqref{eq:discard}. To bound the arithmetic error before setting entries to zero, compare the full product in \eqref{eq:level-update} with $Z_{k+1}^*S_k Z_{k+1}$. We have $\|\widehat Y_{k+1}-Z_{k+1}\|\le q$, $\|\widehat Y_{k+1}\|\le1+q$, and $\|S_k\|<2$. Two applications of the multiplication bound give error less than $12q$ for $q\le1/4$, including replacement of $\widehat Y_{k+1}$ by $Z_{k+1}$. We use the looser bound $16q$. Setting entries to zero can only decrease the Frobenius norm of the arithmetic error, so its remaining spectral norm is at most $16\sqrt n q\le d/8$. The exact block diagonal transformation preserves the earlier zeros. Therefore
\begin{equation}\label{eq:one-level}
 S_{k+1}=Z_{k+1}^*S_k Z_{k+1}+E_{k+1},
 \qquad\|E_{k+1}\|\le9d/16+d/8<d.
\end{equation}
Conjugating by $W_{k+1}$ adds this error in the original coordinates, proving the next invariant.

To see which half-cell contains each child spectrum, use the unitary $Z_0$ from \eqref{eq:direct-rotation} for its parent $C$. The computed child differs from the corresponding exact block of $Z_0^*CZ_0$ by less than $d$, by the same error bounds. The exact child has condition number at most $nK$ by \cref{lem:block}. Bauer--Fike therefore places each computed child eigenvalue within $nKd\le\eta/16$ of the correct exact child spectrum. These discs are disjoint because the parent gap is at least $g_0/2$ and $\eta\le h/(16n)\le g_0/(64n)$. A homotopy along the perturbation gives the matching with multiplicities.

The match through the parent places each child eigenvalue within $5\eta/16$ of its initial label. Equation~\eqref{eq:global-condition} at level $k+1$ places it within $\eta/4$ of the label from the full matrix. If these labels differed, two initial eigenvalues would be within $9\eta/16<g_0$ of one another. Thus the labels agree, and the displacements are too small to cross a dividing line. For a zero or full rank, the entire spectrum belongs to one half-cell. Scalar blocks retain their cells and labels.

A nonscalar block surviving level $D$ would have a cell of side $h\le g_0/4$ and diameter $\sqrt2h<g_0$, containing two distinct initial eigenvalues separated by at least $g_0$. This is impossible. Early termination only replaces $D$ by a smaller completed level in the bounds.

At the terminating level $k\le D$, $Q=\widehat Q_k$ and $T=S_k$. At each multiplication we incur the rounding error of the product as well as the error in the stored factor. Accumulating both gives
\[
 1+\|\widehat Q_j-W_j\|\le(1+q)^{2j}\quad(0\le j\le k),\qquad
 \|Q-W_k\|\le4Dq\le\delta/256.
\]
Here $Dq\le\delta/1024<1/4$. Thus
\begin{equation}\label{eq:Q-close}
 \|Q^*Q-I\|\le2(\delta/256)+(\delta/256)^2<\delta.
\end{equation}
Since $\|T\|<2<3$ and $\|F_k\|\le\mu\le\delta/12$,
\[
 \|A-QTQ^*\|
 \le\delta/3+\delta/12+3(\delta/256)(2+\delta/256)<\delta.
\]
These statements concern the exact matrices represented by the stored entries.
\end{proof}

\section{Precision, probability, and work}\label{sec:completion}
\subsection{Verification of the local error conditions}
To justify \eqref{eq:M0}, we bound the error in forming the regularised matrix. The choice $\zeta\le\nu/(8\gamma\sigma)$ bounds the sampling contribution by $\gamma\sqrt n\zeta\le\nu/8$, and the formation row of \cref{tab:budgets} bounds the elementary conversion and scaling errors by $32\sigma(1+\gamma W)q\le\nu/2$. The combined error from input rounding, Gaussian sampling, scaling by $1/\sqrt n$ and addition is at most $\nu$.

For a shifted sign problem, the initial rounding and each Newton step contribute at most $e=q(KB)^3$, by the kernel bounds with $C_0=2^{20}$ used in \cref{cor:sign}. The computed Newton iterates $Y_j$ satisfy $\|Y_j\|,\|Y_j^{-1}\|\le2KB$ while in the $\rho$-neighbourhood of $X_j$. Thus $q\le1/[16(KB)^2]$ ensures that the inverse bound applies. To keep each iterate in its neighbourhood, we check the stability condition for the induction. With $L=2^{12}K^2B^6$ and $\rho=(4KB)^{-1}$, the stability row gives
\[
 e\le\frac1{2^{31}(N+1)K^5B^{13}}
 \le\frac{\rho}{4(N+1)L^2}.
\]
The accuracy row gives $(N+1)Le\le\tau/16$. The fixed Newton count in \eqref{eq:audit-sign} bounds exact-arithmetic sign truncation by $\tau/4$, as shown in \cref{cor:sign}. The sign error is consequently less than $\tau/2$. Formation of the projector adds at most $q(K+1)\le\tau/4$, proving \eqref{eq:projector-accuracy}. The range-product row similarly gives \eqref{eq:range-error}, and the remaining rows give the trace, transformation and basis bounds already used.

\subsection{A joint probability bound}
Combining the failure probabilities from regularisation, sampling, grid selection and range extraction gives the bound
\begin{equation}\label{eq:total-failure}
 \Pi=\omega+\frac{9n^3}{\gamma^2t^2}
       +\frac{576t^8r^6n^8}{\gamma^8}
       +6n^3a+\frac{4n\eta}{h}+\frac{4n}{M_G}+n^3s^2.
\end{equation}
For the conditioning and gap failures within the norm cutoff, use \eqref{eq:BGVKS}. Exceeding the cutoff on the retained event contributes $\omega$. We count sampler failure only once, even when retaining the samples also bounds the perturbation norm. Conditional on the history up to a range sample request, the new continuous Gaussian sample has its original independent law. Apply \cref{lem:range} before conditioning on acceptance of that sample. There are at most $n-1$ nontrivial splits and each rank is at most $n$, so the total range-failure probability is at most $n^3s^2$. The grid offsets are independent of the initial matrix, so we apply \cref{lem:grid} conditionally on that matrix. The union bound gives total failure probability at most $\Pi$.

For \eqref{eq:explicit-confidence}, take the regularisation norm cutoff $M_{\rm cut}=V_{\rm G}$ and $\omega=0$. The conditioning, gap, sampling, grid and range terms are each at most $\beta=p/8$, so
$
 \Pi\le5\beta=5p/8<p.
$
Outside these failure events, \cref{prop:invariant} gives both target inequalities jointly, and the arithmetic checks in \cref{sec:algorithm} exclude exceptional returns. More generally, any schedule satisfying the deterministic inequalities and $\Pi\le p$ is admissible, with room to vary the individual budgets.

\subsection{Precision and deterministic work}
Let $\Lambda=\log_2(n/(\delta p))$. In \eqref{eq:explicit-confidence}, $\ell=O(\Lambda)$, $W=O(\sqrt\Lambda)$, and $\gamma^{-1},t,r^{-1}$ are bounded by fixed polynomials in $n,\delta^{-1},p^{-1}$. The same is true of every inverse accuracy or margin in \eqref{eq:grid-parameters}--\eqref{eq:local-budgets}. The factors $D$ and $W$ can also be bounded by such polynomials. Thus
\[
 D,J,m_U,N=O(\Lambda),\qquad b_{\rm dir}=O(\Lambda).
\]
All constants are universal, since $C_0$ is a universal kernel constant. Parameter selection uses integer comparisons on $O(\Lambda)$-bit quantities and conservative dyadic rounding. Only scalar operations are needed here. The entry cap in \eqref{eq:R} also has $\log R=O(\Lambda)$.

At each level, current block orders satisfy $\sum m=n$, so $\sum m^3\le n^3$. Sign work is $O(n^3N)$ per level, and full matrix transforms, basis accumulation and QR add $O(n^3)$ per level. There are at most $D$ levels on every execution. Range sampling uses at most $n^2$ entries per level, and regularisation uses $n^2$ entries. The sampler and bit-generation cost is therefore at most $O(n^2(D+1)(J+m_U))$. The one-time scalar constants cost $O(J)$, and range reduction is included in the sampler bound. The two grid offsets use $2\log_2M_G=O(\Lambda)$ additional bits, absorbed in the same bound. Combining all numerical work gives
\begin{align*}
 O\bigl(n^3D(N+1)+n^2(D+1)(J+m_U)\bigr)=O\bigl((n^3+n^2)\Lambda^2\bigr)
 =O(n^3\Lambda^2).
\end{align*}
A rounded real primitive and the production of one unbiased random bit are counted as unit work. The prescribed caps and absence of retries make this an operation bound on every execution, including unsuccessful ones. Together with the scalar case in \cref{sec:random}, this handles every $n\ge1$. Rounding the original requests down only strengthens the guarantees. This proves \cref{thm:main}.

\section{Numerical experiments}\label{sec:numerics}
We measure how the arithmetic work varies with dimension, tolerance and requested failure probability. We use the explicit schedule \eqref{eq:explicit-confidence}, with 32 guard bits added to the sufficient mantissa length in \eqref{eq:precision-choice}. The precision, sampling lengths, Newton counts and subdivision caps are fixed before each run. We stop when all blocks are scalar and never retry a failed run.

\paragraph{Test matrices}
We first vary the dimension over $n=2,4,8,16,32,64$, with $\delta=p=2^{-8}$ and three seeds for each of two matrix families. The generic matrices have real and imaginary integer components drawn uniformly from $\{-4,\ldots,4\}$ and divided by $8n$. To obtain exactly defective matrices, we use $UJ_nU^*$, where $J_n$ is the nilpotent Jordan block and $U=D_\phi H_n/\sqrt n$. Here $H_n$ is the Walsh--Hadamard matrix, and $D_\phi$ is a seeded diagonal of fourth roots of unity. All entries are formed exactly as dyadics, so the defective matrices retain a single Jordan block and have norm one. Both families satisfy $\|A\|\le1$ by construction.

To examine the dependence on tolerance and confidence, we use normal matrices of order four whose eigenvalues lie within $\delta/8$ of $3/8+3i/16$, with the same exact unitary construction. We vary $\delta=2^{-8},2^{-16},2^{-32},2^{-64}$ at $p=2^{-8}$, and then $p=2^{-8},2^{-16},2^{-32}$ at $\delta=2^{-8}$. The eigenvalue cluster shrinks with $\delta$ in the tolerance experiment. In the confidence experiment, each seed's exact input remains fixed. Across all three experiments, reusing the common baseline gives 54 runs in total.

\paragraph{Accuracy}
For each stored output, we check exact upper triangularity and bound the normalised error
\[
 E_F=\delta^{-1}\max\{\|A-QTQ^*\|_F,\|Q^*Q-I\|_F\}.
\]
We evaluate the squared Frobenius norms exactly as rationals. Since the Frobenius norm bounds the spectral norm, $E_F\le1$ certifies both estimates in \cref{thm:main} against the original input. All 54 runs completed without abort and passed these checks, with $\max E_F<0.007608$. The prescribed precisions range from 1695 to 4694 bits. With $\Lambda=\log_2(n/(\delta p))$, they satisfy $63.3<b/\Lambda<111.2$. These are sufficient precisions from the schedule, and the runs leave open how far the precision can be reduced.

\begin{figure}[!t]
\centering
\includegraphics[width=\linewidth]{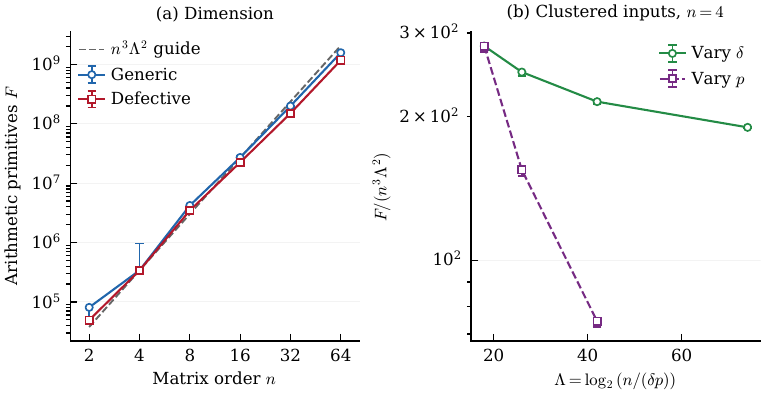}
\caption{Arithmetic work for the explicit sufficient schedule. (a) Primitive counts for generic and exactly defective matrices at $\delta=p=2^{-8}$, with an $n^3\Lambda^2$ guide normalised at $n=4$. (b) Counts divided by $n^3\Lambda^2$ for the tolerance and confidence experiments on clustered normal matrices of order four, where $\Lambda=\log_2(n/(\delta p))$. Points give medians over seeds 40--42, and bars show the minimum--maximum range. The two curves in (b) share their baseline.}
\label{fig:numerics}
\end{figure}

\paragraph{Arithmetic work}
We measure arithmetic work by the numerical primitive count $F$, which includes finite sampling, inversions, QR and full basis updates. Exact sign changes and power-of-two scalings each count as one primitive. We count conversions, comparisons, random-bit draws, exact parameter setup and verification separately. In every run, the random-bit count is less than $5.1\%$ of $F$.

Figure~\ref{fig:numerics}(a) compares the generic and defective families with an $n^3\Lambda^2$ guide, normalised at $n=4$. At $n=64$, the median primitive counts are $1.58\times10^9$ for generic matrices and $1.18\times10^9$ for defective matrices. Panel (b) shows $F/(n^3\Lambda^2)$ for the tolerance and confidence experiments. Over the tested ranges, the unscaled median count grows by factors of $11.41$ and $1.44$ as $\delta$ and $p$ decrease, respectively, while the normalised counts in both experiments fall.

\paragraph{Implementation}
A compiled C\# backend rounds each real operation to nearest, with ties to even and arbitrary-size exponents. Seeds 40--42 specify reproducible keyed pseudorandom bit streams that simulate the independent-bit model of the analysis. We used a host with two 16-core Intel Xeon Gold 6246R processors at nominal 3.4\,GHz and 376\,GiB RAM, with Python 3.10.12 and PowerShell 7.6.9/.NET 10.0.14. Eight~independent workers each used one pinned core. With these workers running concurrently, the 54 runs took 62.7\,minutes, including exact verification.

\section*{AI declaration}
As documented in the repository, the solution to this problem was found by ChatGPT. The author wrote the manuscript and explanations, and then used ChatGPT for minor language editing only. The author also used ChatGPT to help run the code for the numerical experiment, the author carefully checked all AI output. The author takes responsibility for all content.

\bibliographystyle{siamplain}
\begingroup
\interlinepenalty=10000
\bibliography{references}
\endgroup
\end{document}